%% file: paper2.tex
\documentclass[11pt]{article}
\usepackage[a4paper,margin=25mm]{geometry}
\usepackage{amsmath,amssymb,bm}
\usepackage{booktabs}
\usepackage{graphicx}
\usepackage{subcaption}
\usepackage{siunitx}
\usepackage{xcolor}
\usepackage{hyperref}
\usepackage{microtype}
\usepackage[numbers,sort&compress]{natbib}
\graphicspath{{figures/}}
\newcommand{\Linfty}{L^{\infty}}
\newcommand{\vect}[1]{\bm{#1}}
\title{Central-Hermite Sensing and Collision for Frame-Robust Order-Resolved Relaxation on D3Q125}
\author{Bj\o rn Wu\\
\small Independent Researcher, Oslo, Norway}
\date{}
\begin{document}
\maketitle
\input{sections/abstract}
\input{sections/introduction}
\input{sections/formulation}
\input{sections/models}
\input{sections/homogeneous}
\input{sections/transport_collision}
\input{sections/matrix}
\input{sections/long_time}
\input{sections/transport_limitations}
\input{sections/discussion}
\input{sections/conclusions}
\bibliographystyle{unsrtnat}
\bibliography{references}
\end{document}

%% file: sections/abstract.tex
\begin{abstract}
Raw-Hermite sensing and collision on a fixed discrete-velocity set can convert a uniform translation into artificial coupling between nominally distinct nonequilibrium orders. We develop a central-Hermite formulation for a D3Q125 kinetic model with order-resolved log-Gaussian relaxation and compare three variants: raw sensing/raw collision (A), central sensing/raw collision (B), and central sensing/central collision (C). In homogeneous translated second-order perturbations, model C preserves third- and fourth-order modal purity to machine precision, whereas A and B develop boost-dependent cross-order content. Across a grid--CFL--boost matrix, model C reduces the post-transport collision frame discrepancy relative to A by \(65.342\)--\(98.102\%\) (median \(81.131\%\)) in the total relative \(L^\infty\) measure. Long-time calculations remain positive and conservative to numerical precision, although the accumulated benefit is configuration dependent because transport continually re-injects frame error. A transport study further reveals a clear trade-off: central-Hermite interface reconstruction strongly suppresses the third-order discrepancy but amplifies the fourth-order discrepancy. The fully central-Hermite collision therefore substantially reduces collision-induced cross-order frame discrepancy, while residual dependence remains due to discrete transport and finite velocity-space representation. This moment-space improvement does not by itself establish a comparable reduction in macroscopic Galilean transport error.
\end{abstract}

%% file: sections/introduction.tex
\section{Introduction}

Raw Hermite coefficients are sensitive to uniform translation. A physical second-order perturbation can therefore acquire apparent third- and fourth-order content after a boost. When raw coefficients drive both an adaptive sensor and an order-dependent collision operator, this coordinate mixing changes the numerical relaxation even though the underlying physical perturbation is unchanged.

Frame dependence in finite discrete-velocity models is closely connected to the accuracy with which the velocity set represents higher-order moments. Classical analyses of lattice-Boltzmann dispersion, isotropy, stability, and Galilean invariance identify velocity-dependent transport errors when the retained moment content is insufficient~\cite{lallemand2000theory,nie2008galilean}. Systematic Hermite projection combined with sufficiently accurate Gauss--Hermite quadrature provides a route to higher-order discrete kinetic models~\cite{shan2006kinetic}, while regularized formulations control unsupported nonequilibrium content before collision~\cite{latt2006regularized,zhang2006efficient}.

A complementary line of development performs collision in moments shifted by the local hydrodynamic velocity. Cascaded or central-moment lattice-Boltzmann methods use this local frame to reduce translation-induced artifacts and improve robustness~\cite{geier2006cascaded}. Subsequent work has clarified the role of higher-order Hermite content in central-moment formulations and the relation between the continuous Maxwellian attractor and its discrete representation~\cite{derosis2019role}.

This work extends two earlier developments: A log-Gaussian scale-space limiter was first introduced as a smooth probability partition between continuum and ballistic transport in logarithmic Knudsen space~\cite{wu2026loggaussian}. That work established the complementary error-function weighting used here and provided DVM/BGK-informed calibration evidence for broad transition parameters. A subsequent fixed-D3Q125 study introduced hierarchical, order-resolved log-Gaussian relaxation, in which second-, third-, and fourth-order nonequilibrium sectors activate distinct relaxation spectra while sharing the same discrete-velocity representation~\cite{wu2026hierarchical}. The present work retains that order-resolved relaxation architecture but changes the sensing and collision basis from raw Hermite to central Hermite in order to isolate and suppress translation-induced cross-order coupling.

Three otherwise identical formulations are compared: model A uses raw sensing and raw collision; model B uses central sensing and raw collision; model C uses central sensing and central collision. The term frame-robust is used rather than Galilean invariant, because both the finite velocity set and the transport discretization remain tied to the laboratory frame.

The main contributions are: (i) a central-Hermite sensor resolved by nonequilibrium order; (ii) a central-Hermite collision map that prevents translation-induced cross-order relaxation in the homogeneous modal-purity test; (iii) controlled A/B/C comparisons that isolate sensing from collision; and (iv) a transport study revealing a pronounced third-order/fourth-order trade-off under central-Hermite interface reconstruction.

%% file: sections/formulation.tex
\section{Central-Hermite formulation}
The kinetic state is represented by populations $f_i(\vect{x},t)$ associated with the fixed D3Q125 velocity set $\vect{\xi}_i$. Macroscopic density, velocity, and temperature are
\begin{equation}
\rho=\sum_i f_i,\qquad
\rho\vect{u}=\sum_i f_i\vect{\xi}_i,\qquad
T=\frac{1}{3\rho}\operatorname{tr}\vect{C}^{(2)},
\end{equation}
where $\vect{c}_i=\vect{\xi}_i-\vect{u}$ and
\begin{equation}
C^{(n)}_{\alpha_1\cdots\alpha_n}
=\sum_i f_i c_{i\alpha_1}\cdots c_{i\alpha_n},\qquad n=2,3,4.
\end{equation}
With $I_{\alpha\beta}=\delta_{\alpha\beta}$, the nonequilibrium central-Hermite channels used in the implementation are
\begin{align}
\Delta^{(2)}_{\alpha\beta}
&=C^{(2)}_{\alpha\beta}-\rho T I_{\alpha\beta},\\
\Delta^{(3)}_{\alpha\beta\gamma}
&=C^{(3)}_{\alpha\beta\gamma},\\
\Delta^{(4)}_{\alpha\beta\gamma\delta}
&=C^{(4)}_{\alpha\beta\gamma\delta}
-T\,\mathcal{S}_{6}\!\left(I_{\alpha\beta}C^{(2)}_{\gamma\delta}\right)
+\rho T^2\,\mathcal{S}_{3}\!\left(I_{\alpha\beta}I_{\gamma\delta}\right).
\end{align}
Here $\mathcal{S}_{6}$ denotes the six distinct pairings of one identity tensor with $C^{(2)}$, and $\mathcal{S}_{3}$ denotes the three distinct pairings of two identity tensors. 

Here $\mathcal S_6(I\otimes C^{(2)})$ denotes the sum over the six distinct placements of the two indices of $I$ and the two indices of $C^{(2)}$.

For example,
\[
\mathcal S_3(\boldsymbol u\otimes C^{(2)})_{\alpha\beta\gamma}
=
u_\alpha C^{(2)}_{\beta\gamma}
+
u_\beta C^{(2)}_{\alpha\gamma}
+
u_\gamma C^{(2)}_{\alpha\beta}.
\]
Thus $\Delta^{(4)}$ is the fourth central-Hermite coefficient rather than the ordinary fourth central moment. At local equilibrium, $C^{(2)}=\rho T I$ and $C^{(4)}=\rho T^2\mathcal S_3(I\otimes I)$; substitution into the definition above gives $\Delta^{(4)}=0$.

The conversion used for population reconstruction is written most compactly in terms of raw moments $M^{(n)}=\sum_i f_i\vect{\xi}_i^{\otimes n}$:
\begin{align}
\vect{M}^{(2)}&=\vect{C}^{(2)}+\rho\vect{u}^{\otimes2},\\
\vect{M}^{(3)}&=\vect{C}^{(3)}+\mathcal{S}_{3}(\vect{u}\otimes\vect{C}^{(2)})+\rho\vect{u}^{\otimes3},\\
\vect{M}^{(4)}&=\vect{C}^{(4)}+\mathcal{S}_{4}(\vect{u}\otimes\vect{C}^{(3)})
+\mathcal{S}_{6}(\vect{u}^{\otimes2}\otimes\vect{C}^{(2)})+\rho\vect{u}^{\otimes4}.
\end{align}
The raw Hermite coefficients are then
\begin{align}
\vect{a}^{(2)}&=\vect{M}^{(2)}-\rho\vect{I},\\
\vect{a}^{(3)}&=\vect{M}^{(3)}-\mathcal{S}_{3}(\vect{I}\otimes\rho\vect{u}),\\
\vect{a}^{(4)}&=\vect{M}^{(4)}-\mathcal{S}_{6}(\vect{I}\otimes\vect{M}^{(2)})
+\rho\mathcal{S}_{3}(\vect{I}\otimes\vect{I}).
\end{align}
These triangular identities make explicit why a uniform translation mixes a physical second-order perturbation into raw third- and fourth-order coefficients, whereas the local central channels are unchanged by that translation.

The D3Q125 velocity set is the tensor product of a five-node one-dimensional Gauss--Hermite quadrature,
\[
\mathcal V_{125}=\mathcal V_5\otimes\mathcal V_5\otimes\mathcal V_5,
\]
with tensor-product weights \(w_i\). The one-dimensional rule is exact for polynomials through degree nine with respect to the standard-Gaussian weight. The local equilibrium used in the computations is the fourth-order discrete Hermite projection of the Maxwellian associated with \((\rho,\boldsymbol u,T)\), rather than a continuously sampled Maxwellian.

%% file: sections/models.tex
\section{Compared formulations}
The three formulations are
\begin{align}
\text{A}:&\quad \text{raw sensor}+\text{raw collision},\\
\text{B}:&\quad \text{central sensor}+\text{raw collision},\\
\text{C}:&\quad \text{central sensor}+\text{central collision}.
\end{align}
All variants use the same velocity set, equilibrium, transport operator, time step, boundaries, and order-dependent relaxation spectrum. Consequently, A versus B isolates the effect of sensing, while B versus C isolates the collision basis.

For the central sensor, the dimensionless order contributions are
\begin{align}
\chi_2&=\frac{\|\Delta^{(2)}\|_F}{\max(\rho T,\varepsilon)},\\
\chi_3&=\frac{\|\Delta^{(3)}\|_F}{\max(\rho T,\varepsilon)\sqrt{\max(T,\varepsilon)}},\\
\chi_4&=\frac{\|\Delta^{(4)}\|_F}{\max(\rho T,\varepsilon)\max(T,\varepsilon)},\\
\chi&=\chi_2+\chi_3+\chi_4,
\end{align}
with Frobenius norms, unit coefficients for all three orders, and $\varepsilon=10^{-14}$. The raw variant uses the corresponding raw-Hermite deviations with the same normalization. Central sensing therefore removes the direct algebraic boost dependence of the local nonequilibrium indicator; it does not remove frame dependence introduced later by the finite velocity set or transport discretization.

The order-resolved effective indicator supplied to the relaxation spectrum is
\[
K_n=\left(K_\rho^p+K_T^p+K_u^p+\chi_n^p\right)^{1/p},
\qquad n=2,3,4,
\]
with \(p=8\). For the periodic one-dimensional tests,
\[
K_\rho=\lambda\frac{|\partial_x\rho|}{\max(\rho,\varepsilon)},\quad
K_T=\lambda\frac{|\partial_xT|}{\max(T,\varepsilon)},\quad
K_u=\lambda\frac{\\|\partial_x\boldsymbol u\\|_2}{\sqrt{\max(T,\varepsilon)}}.
\]
The parameter \(\lambda\) is the prescribed gradient-sensor length scale. In the TNE-only limit \(\lambda=0\), the effective indicator is driven solely by \(\chi_n\).

For each order $n\in\{2,3,4\}$, the continuum weight is the log-Gaussian transition
\begin{equation}
w_{c,n}(K_n)=\frac{1}{2}\operatorname{erfc}\!\left[
\frac{\log(\max(K_n,10^{-14})/K_{0,n})}{\sqrt{2}\,\sigma_n}
\right],\qquad w_{k,n}=1-w_{c,n}.
\end{equation}
The calibrated reference-step relaxation factor is
\begin{equation}
s_n^{\mathrm{ref}}(K_n)=w_{c,n}s_{c,n}+w_{k,n}s_{k,n},
\end{equation}
with
\begin{center}
\begin{tabular}{c c c c c}
\toprule
$n$ & $K_{0,n}$ & $\sigma_n$ & $s_{c,n}$ & $s_{k,n}$\\
\midrule
2 & 0.050 & 2.0 & 1.00 & 0.20\\
3 & 0.030 & 2.5 & 1.00 & 0.10\\
4 & 0.015 & 3.0 & 1.00 & 0.05\\
\bottomrule
\end{tabular}
\end{center}
For a time step $\Delta t$, the implementation preserves the reference survival fraction:
\begin{equation}
s_n(\Delta t)=1-\left(1-s_n^{\mathrm{ref}}\right)^{\Delta t/\Delta t_{\mathrm{ref}}},
\qquad \Delta t_{\mathrm{ref}}=1.0938161705673147\times10^{-3},
\end{equation}
with the limiting value $s_n=1$ retained exactly when $s_n^{\mathrm{ref}}=1$. The configured factors lie in $[0,1]$; the central collision routine itself validates the more general interval $0\le s_n\le2$.

The fully central collision updates only the nonconserved channels,
\begin{equation}
\Delta^{(n),+}=(1-s_n)\Delta^{(n)},\qquad n=2,3,4.
\end{equation}
Density, velocity, and the trace-defined temperature are held fixed. The post-collision central moments are reconstructed from the relaxed channels, converted to raw moments by the identities above, converted to raw Hermite coefficients, and finally mapped back to populations. This construction preserves mass and momentum by design. Moreover,
\[
\operatorname{tr}\Delta^{(2)}
=
\operatorname{tr}C^{(2)}-3\rho T
=
0
\]
because \(T=\operatorname{tr}C^{(2)}/(3\rho)\). Multiplication by \(1-s_2\) therefore leaves the second-moment trace unchanged, so the trace-defined temperature and local discrete kinetic energy are preserved to floating-point accuracy.

%% file: sections/homogeneous.tex
\section{Homogeneous translated perturbations}
For a translated second-order perturbation, the central total sensor remains $0.0244948974$ to numerical tolerance, while the third- and fourth-order relaxation factors remain equal to one. Model C preserves modal purity: through boost $U=0.4$, the post-collision central third- and fourth-order residuals stay at round-off level. At $U=0.4$, model B produces residuals of approximately $7.48\times10^{-3}$ and $5.77\times10^{-3}$, respectively. The post-collision second-order amplitude remains $7.066848\times10^{-3}$ for every boost.

\begin{figure}[t]
\centering
\includegraphics[width=0.95\linewidth]{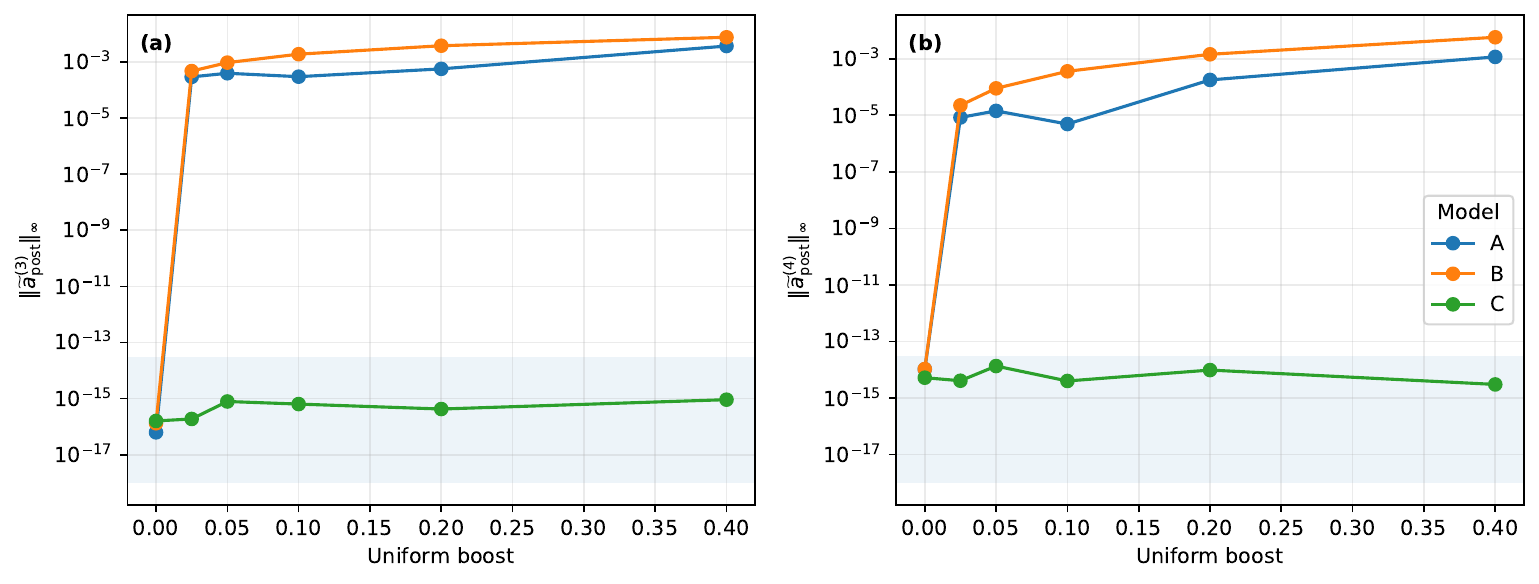}
\caption{Post-collision central third- and fourth-order residuals under uniform translation. Model C remains at the numerical round-off floor, whereas models A and B develop boost-dependent cross-order content. The zero-boost values are included as round-off references.}
\label{fig:homogeneous-purity}
\end{figure}

%% file: sections/transport_collision.tex
\section{Transport-generated nonequilibrium}
The diagnostic uses a periodic one-dimensional smooth compression wave on a domain of length $L=1$ with $N=24$ cells. The initial state is the corresponding fourth-order discrete Hermite equilibrium with
\begin{align}
\rho(x)&=1+0.08\cos(2\pi x/L),\\
u_x(x)&=U+0.06\sin(2\pi x/L),\qquad u_y=u_z=0,\\
T(x)&=1+0.08\cos(2\pi x/L+\pi/4).
\end{align}
The macroscopic-gradient sensor scale is $\lambda=0.01$, the boundary condition is periodic, and first-order population upwinding is used with CFL number $0.4$. The time step is therefore
\begin{equation}
\Delta t=0.4\,\frac{\Delta x}{\max_i|\xi_{ix}|},\qquad \Delta x=L/N.
\end{equation}

Because the initial state is locally equilibrated, the informative sequence is one transport step followed by one collision step.

Here $\lambda$ multiplies the normalized macroscopic gradients in the adaptive sensor and is not an independently prescribed molecular collision time.

Let $\chi_U(x)$ denote the central total-TNE profile at uniform boost $U$, and let $\chi_0(x)$ be the zero-boost profile at the same stage. The reported direct relative discrepancy is
\begin{equation}
D_{\infty}(U)=
\frac{\|\chi_U-\chi_0\|_{\infty}}
{\max(\|\chi_0\|_{\infty},10^{-14})}.
\end{equation}

This quantity is a central-TNE moment-space diagnostic; it is not a direct measure of boost-dependent error in a macroscopic transport coefficient, phase speed, or decay rate.
Transport establishes a baseline discrepancy; raw collision strongly amplifies it; central collision suppresses that amplification. After transport and collision, $D_{\infty}$ at $U=0.1$ is approximately 1.295 for A, 0.2765 for B, and 0.04384 for C. At $U=0.2$, the corresponding values are approximately 4.139, 0.5230, and 0.07858. Thus, the fully central collision acts as a strong local correction after transport has generated nonequilibrium.

%% file: sections/matrix.tex
\section{Grid--CFL--boost matrix}
For the total relative $\Linfty$ central-TNE frame discrepancy, model C improves on A in every nonzero-boost configuration tested. These reductions quantify a moment-space collision diagnostic and should not be read as direct reductions of macroscopic Galilean transport error.

The reduction ranges from 65.342\% to 98.102\%, with median 81.131\%. The second-order channel is unchanged to numerical precision and is therefore omitted from Table~1 rather than being a missing entry. The third-order reduction ranges from 84.535\% to 93.604\% (median 88.790\%), and the fourth-order reduction ranges from 42.607\% to 96.526\% (median 67.321\%).

\begin{table}[t]
\centering
\caption{Reduction of model C relative to model A in the one-step grid--CFL--boost matrix.}
\label{tab:matrix-summary}
\begin{tabular}{lrrr}
\toprule
Metric & Minimum (\%) & Median (\%) & Maximum (\%)\\
\midrule
Total relative $\Linfty$ & 65.342 & 81.131 & 98.102\\
Third order & 84.535 & 88.790 & 93.604\\
Fourth order & 42.607 & 67.321 & 96.526\\
\bottomrule
\end{tabular}
\end{table}

\begin{figure}[t]
\centering
\includegraphics[width=0.92\linewidth]{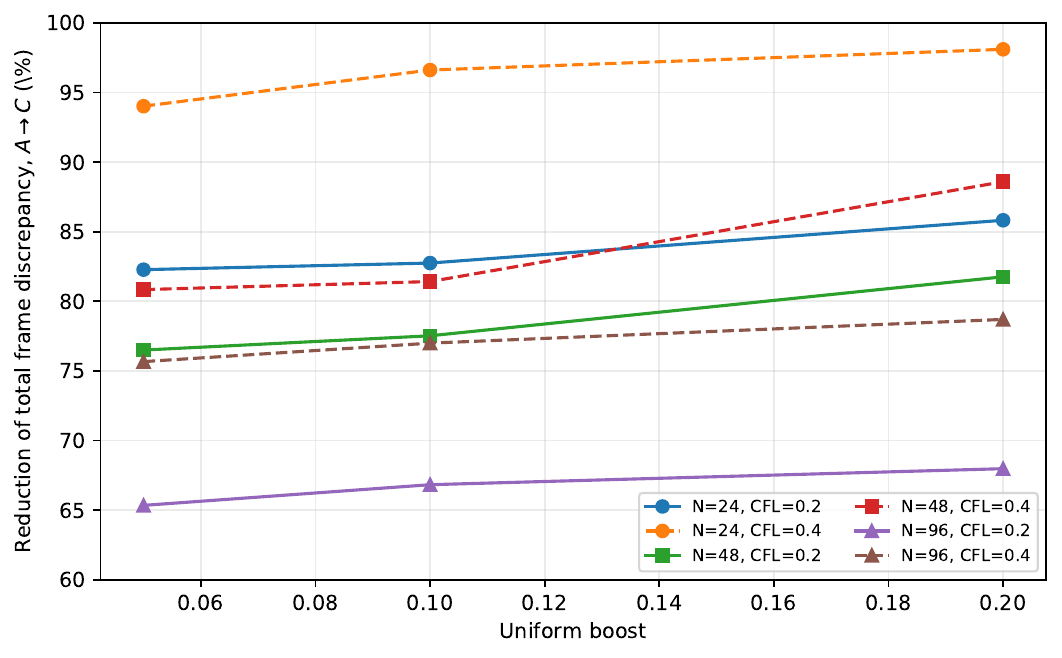}
\caption{Reduction of the total frame discrepancy from model A to model C across grid resolution, CFL number, and uniform boost. The vertical axis is restricted to 60--100\% to resolve configuration dependence.}
\label{fig:matrix-reduction}
\end{figure}

%% file: sections/long_time.tex
\section{Long-time behaviour}
The minimum population over the long-time matrix is $3.617\times10^{-7}$. The maximum relative mass error is $4.012\times10^{-14}$, the maximum absolute momentum error is $5.613\times10^{-13}$, and the maximum relative energy error is $2.758\times10^{-14}$. The final-time reduction of C relative to A ranges from $-0.324\%$ to 19.013\% (mean 7.332\%). The time-integrated reduction ranges from 0.343\% to 40.118\% (mean 16.737\%). The weaker long-time benefit is consistent with transport repeatedly generating new frame error and underscores that collision-level correction need not control the dominant macroscopic transport error of the full discretization.

\begin{figure}[t]
\centering
\includegraphics[width=0.90\linewidth]{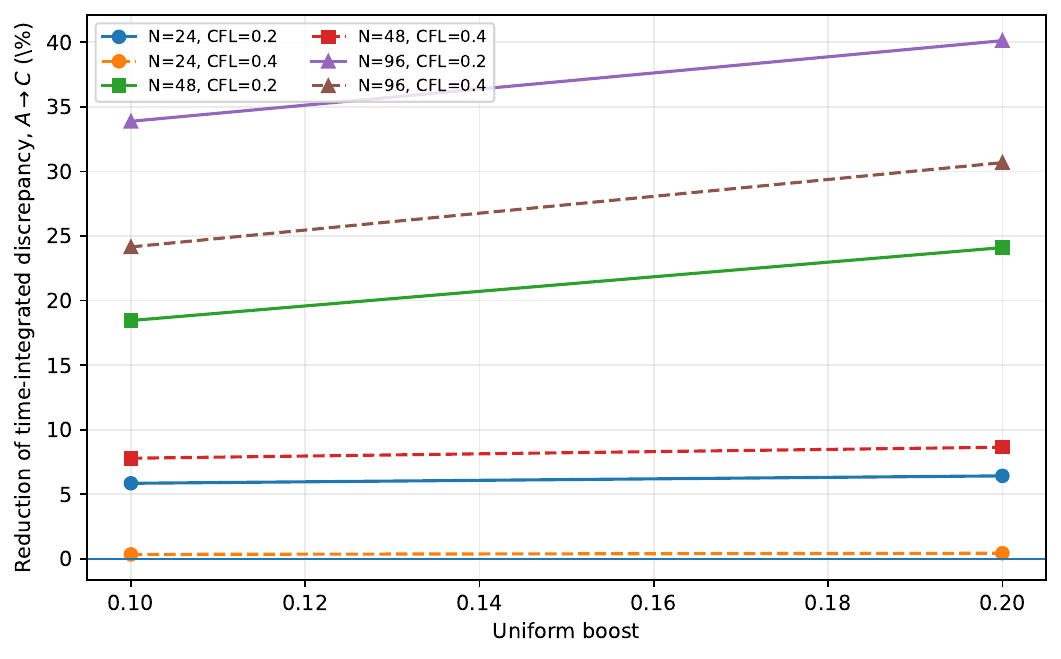}
\caption{Reduction of time-integrated frame discrepancy for model C relative to model A. The benefit is positive throughout the tested matrix but depends on resolution, CFL number, and boost.}
\label{fig:long-time}
\end{figure}

The slightly negative lower endpoint is a final-time, phase-sensitive value rather than a persistent degradation: the time-integrated reduction remains positive in every tested configuration. This is consistent with a small crossing of the two discrepancy histories near the selected terminal time, while the accumulated discrepancy over the complete trajectory remains lower for model C.

%% file: sections/transport_limitations.tex
\section{Transport limitation study}
Population-wise MUSCL reconstruction, SSP-RK2 time integration, and limiter changes do not systematically reduce frame discrepancy. Central-Hermite interface reconstruction reduces the third-order discrepancy by 87.1--94.1\%, but amplifies the fourth-order discrepancy by approximately 141--588\% in the tested matrix. A shared interface frame does not remove this fourth-order increase.

The maximum absolute round-trip reconstruction errors are $2.165\times10^{-15}$ for order two, $1.685\times10^{-15}$ for order three, and $1.443\times10^{-14}$ for order four. The residual low-order leakage is $5.465\times10^{-15}$, and transporting that residual leaves the result unchanged. These audits exclude a faulty fourth-order translation formula, a mismatched shared frame, or loss of a pre-existing high-order residual as the dominant cause. The observations are consistent with noncommutation between spatial interpolation and the nonlinear central-Hermite map, potentially compounded by finite-quadrature aliasing.

\begin{figure}[t]
\centering
\includegraphics[width=0.83\linewidth]{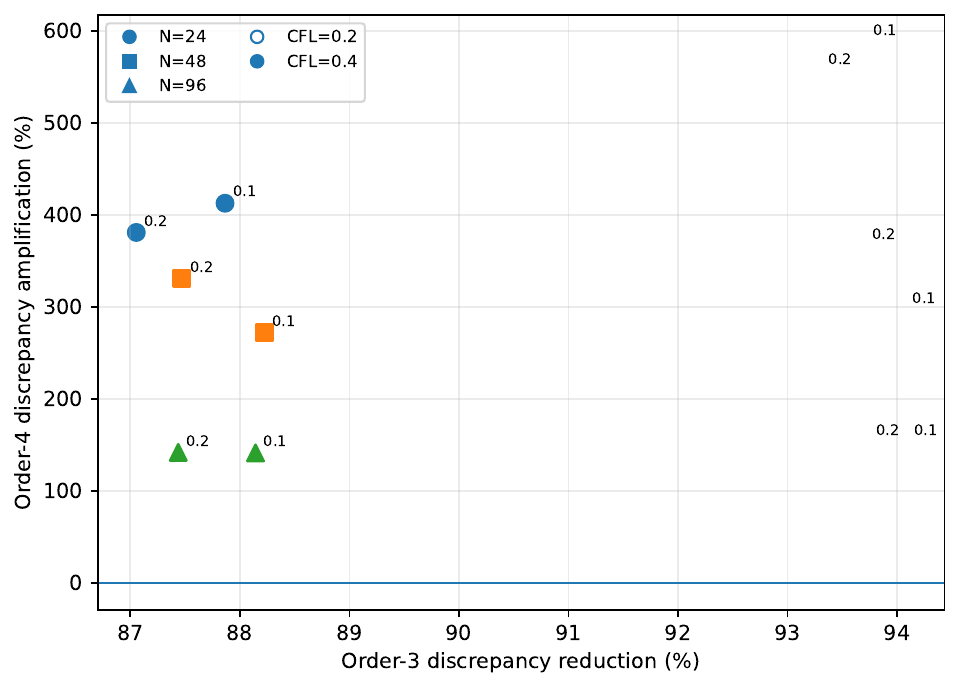}
\caption{Trade-off produced by central-Hermite interface reconstruction. Third-order frame discrepancy is strongly reduced, while fourth-order discrepancy is amplified. Marker shape identifies grid resolution and marker fill distinguishes CFL number; point annotations identify the boost.}
\label{fig:transport-tradeoff}
\end{figure}

%% file: sections/discussion.tex
\section{Discussion}
Central sensing removes the direct algebraic boost dependence from the adaptive indicator, but central collision is required to preserve modal purity. After transport, the fully central collision provides a strong local reduction of translation-induced cross-order coupling. It does not make the complete solver exactly Galilean invariant because the fixed velocity representation and transport discretization remain frame dependent.

The large percentage reductions reported above concern post-transport, post-collision central-TNE discrepancy. They demonstrate that the central collision basis avoids further amplification of translation-induced cross-order content. They do not demonstrate that the same percentage reduction carries over to macroscopic Galilean transport error. Consequently, strong improvement of the collision-space diagnostic does not establish that collision error is the dominant contribution to boost-dependent macroscopic transport error.

The transport study supplies a useful negative result. Higher formal order does not guarantee greater frame robustness, and independent reconstruction of nonlinear central moments can exchange a large third-order improvement for severe fourth-order aliasing. This observation argues against treating moment-wise central reconstruction as a universally benign upgrade.

The scope of the present study is deliberately restricted. The numerical evidence concerns the tested D3Q125 representation, smooth periodic perturbations, grid resolutions, CFL numbers, and uniform-boost range. No external kinetic-reference solution is used here, and the calculations do not yet cover strongly nonequilibrium states, discontinuities, or non-periodic boundary treatments. The central-Hermite implementation is a research prototype rather than an optimized production kernel, so the reported frame-robustness results should not be interpreted as a computational-cost benchmark. These limitations do not affect the controlled A/B/C comparison, but they delimit the range over which the quantitative reductions can presently be claimed.

%% file: sections/conclusions.tex
\section{Conclusions}
The fully central-Hermite formulation provides a substantial and consistent reduction of collision-induced cross-order frame discrepancy across the tested configurations. It prevents translation-induced cross-order relaxation in the homogeneous modal-purity tests and reduces the post-transport collision discrepancy by 65.342--98.102\% relative to the fully raw model across the tested grid--CFL--boost matrix, while maintaining positivity and conservation to numerical precision in the tested cases. This improvement is specific to the reported moment-space and collision-level diagnostics and does not by itself imply a comparable reduction in macroscopic Galilean transport error. Long-time benefits are smaller and configuration dependent because transport continually re-injects frame error. The remaining evidence points primarily to discrete transport and finite velocity-space representation, rather than to the central collision map, as the dominant sources of residual frame dependence in the tested solver.